\documentclass[11pt]{article}
\usepackage[T1]{fontenc}
\usepackage{amsfonts}
\usepackage{mathtools}
\usepackage{color}
\usepackage{amsmath,amssymb,comment}
\usepackage{hyperref}
\hypersetup{%
  colorlinks=true,
  linkcolor=black,
  citecolor=black,
  urlcolor=blue,
  pdftitle={Uniform Lp estimates for solutions to the inhomogeneous 2D Navier--Stokes equations and application to a chemotaxis--fluid system with local sensing},
  pdfauthor={},
  pdfkeywords={},
  bookmarksopen=true,
}

\newtheorem{theo}{Theorem}[section]
\newtheorem{lem}[theo]{Lemma}

\newcommand{\mysection}[1]{\section{#1} \setcounter{equation}{0}}
    \makeatletter
\def\@fnsymbol#1{\ensuremath{\ifcase#1\or *\or \ddagger\or
   \mathsection\or \mathparagraph\or \|\or **\or \dagger\dagger
   \or \ddagger\ddagger \else\@ctrerr\fi}}
    \makeatother
\newcommand{\proof}{{\sc Proof.} \quad}
\newcommand{\proofc}{{\sc Proof} \ }
\newcommand{\be}{\begin{equation} \label}
\newcommand{\ee}{\end{equation}}
\newcommand{\bea}{\begin{eqnarray}\label}
\newcommand{\eea}{\end{eqnarray}}
\newcommand{\bas}{\begin{eqnarray*}}
\newcommand{\eas}{\end{eqnarray*}}
\newcommand{\bit}{\begin{itemize}}
\newcommand{\eit}{\end{itemize}}
\newcommand{\qed}{\hfill$\Box$ \vskip.2cm}
\newcommand{\nn}{\nonumber}
\newcommand{\R}{\mathbb{R}}

\newcommand{\pO}{\partial\Omega}

\newcommand{\eps}{\varepsilon}

\newcommand{\hra}{\hookrightarrow}
\newcommand{\io}{\int_\Omega}
\newcommand{\na}{\nabla}
\newcommand{\Del}{\Delta}

\newcommand{\al}{\alpha}

\newcommand{\pa}{\partial}
\newcommand{\bom}{\overline{\Omega}}
\newcommand{\Om}{\Omega}

\newcommand{\ov}{\overline}

\newcommand{\hs}{\hspace*}

\newcommand{\proj}{{\mathcal{P}}}
\newcommand{\vp}{\varphi}
\newcommand{\lbal}{\left\{ \begin{array}{l}}
\newcommand{\lball}{\left\{ \begin{array}{ll}}
\newcommand{\ear}{\end{array} \right.}

\newcommand{\abs}{\\[5pt]}
\newcommand{\Abs}{\\[5mm]}

\newcommand{\tm}{T_{\mathrm{max}}}

\newcommand{\ow}{\ov{w}}

\newcommand{\onz}{\ov{n}_0}

\newcommand{\ure}{\mathrm e}
\newcommand{\diff}{\mathrm d}
\newcommand{\dx}{\,\diff x}
\newcommand{\ds}{\,\diff s}
\newcommand{\dsigma}{\,\diff \sigma}
\newcommand{\mc}{\mathcal}
\newcommand{\defs}{\coloneqq}
\newcommand{\ddt}{\frac{\diff}{\diff t}}

\begin{document}
\enlargethispage{10mm}
\title{Uniform $L^p$ estimates for solutions to the inhomogeneous 2D Navier--Stokes equations and application to a chemotaxis--fluid system with local sensing}

\author{
{Mario Fuest\footnote{fuest@ifam.uni-hannover.de, corresponding author}}\\
{\small Leibniz Universität Hannover, Institut für Angewandte Mathematik,}\\
{\small 30167 Hannover, Germany}
\and
Michael Winkler\footnote{michael.winkler@math.uni-paderborn.de}\\
{\small Institut f\"ur Mathematik, Universit\"at Paderborn,}\\
{\small 33098 Paderborn, Germany}
}
\date{}
\maketitle
\begin{abstract}
\noindent
  The chemotaxis-Navier--Stokes system
  \begin{equation*}\label{1}
  \left\{
  \begin{array}{rcl}
	n_t+u\cdot\nabla n
	&=& \Delta \big(n c^{-\alpha} \big), \\[1mm]
	c_t+ u\cdot\nabla c
	&=& \Delta c -nc,\\[1mm]
    	u_t + (u\cdot\nabla) u
	&=&\Delta u+\nabla P + n\nabla\Phi,
	\qquad \nabla\cdot u=0,
	\end{array} \right.
  \end{equation*}
  modelling the behavior of aerobic bacteria in a fluid drop,
  is considered in a smoothly bounded domain $\Omega \subset \mathbb R^2$.
  For all $\alpha > 0$ and all sufficiently regular $\Phi$, we construct global classical solutions
  and thereby extend recent results for the fluid-free analogue to the system coupled to a Navier--Stokes system.  \abs
  As a crucial new challenge, our analysis requires a~priori estimates for $u$ at a point in the proof when knowledge about $n$
  is essentially limited to the observation that the mass is conserved.
  To overcome this problem, we also prove new uniform-in-time $L^p$ estimates for solutions to the inhomogeneous 
  Navier--Stokes equations merely depending on the space-time $L^2$ norm of the force term raised to an arbitrary small power.\\[2pt]
\end{abstract}

%
%
%
%
%
%
%
%
%
\section{Introduction}\label{intro}
{\bf The two-dimensional Navier--Stokes equations with a force term in $L^\infty$-$L^1$.} \quad
In the first part of the present paper,
we are concerned with the regularity of solutions $(v, Q)$ of the two-dimensional incompressible Navier--Stokes equations
\be{0v}
\lball
v_t + (v\cdot\na) v = \Del v + \na Q + f(x,t),
\quad \na \cdot v=0,
\qquad & x\in\Om, \ t\in (0,T), \\[1mm]
v=0,
\qquad & x\in\pO, \ t\in (0,T), 
\ear
\ee
where $\Omega \subset \R^2$ is a smoothly bounded domain and $f$ is a given external force.
In contrast to the three-dimensional setting
where strong a~priori estimates remain elusive even if $f \equiv 0$ (cf.\ \cite{LerayMouvementLiquideVisqueux1934}, \cite{sohr}),
the 2D case is quite well understood even for forces $f$ with rather weak integrability and regularity properties.
For instance, if $f \in L^1((0, T); L^2(\Om;\R^2)) + L^2((0, T); W^{-1, 2}(\Om;\R^2))$
then for any sufficiently regular $v_0$ there exists a unique so-called strong solution $v$ of \eqref{0v} with $v(\cdot, 0) = v_0$ inter alia belonging to $L^4(\Omega \times (0, T);\R^2)$ \cite[Theorem~4.2.1]{sohr},
and further well-posedness results for $f \in L^r((0, T); W^{-1, q}(\Om;\R^2))$ have been obtained in \cite{CasasKunischWellposednessEvolutionaryNavierStokes2021} for other choices of $r$ and $q$ as well.
Moreover, \cite[Proposition~1.1 and Proposition~1.3]{win_lp} show that
\be{winprop1}
  \begin{cases}
    \text{if $\frac1r + \frac1q \le 1$ and $q > 1$, bounds for $f$ in $L^r((0, T); L^q(\Om))$ (and suitable ones for $v(\cdot, 0)$)} \\
    \text{imply uniform-in-time $L^p(\Om)$ estimates for $v$ for all $p \in [1, \infty)$}.
  \end{cases}
\ee
Let us compare the latter result to the two-dimensional Stokes equations, that is, \eqref{0v} without the nonlinear convection term $(v \cdot \na) v$.
For $p \in (1, \infty)$, the Stokes operator $A_p$ on $L_\sigma^p(\Om) \defs \{\varphi \in L^p(\Om; \R^2) \mid \nabla \cdot \varphi = 0 \text{ in } \mc D'(\Om) \}$
generates a bounded analytic semigroup \cite{giga1981_the_other}.
Therefore, the analogue of \eqref{winprop1} for the Stokes equations follows by well-known semigroup estimates,
while using somewhat more subtle arguments (cf.\ \cite[Lemma 2.5]{WangXiangGlobalExistenceBoundedness2015} and also Lemma~\ref{lem2} below, for instance)
one can show that \eqref{winprop1} continues to hold in the borderline case $r=\infty$ and $q=1$, i.e., if
\be{1.1}
\sup_{t\in (0,T)} \|f(\cdot,t)\|_{L^1(\Om)} <\infty.
\ee
However, whether or not \eqref{1.1} is sufficient for uniform $L^p$ bounds also for solutions to the full Navier--Stokes equations seems to be an open problem.
The key challenge for deriving such a result lies in the fact that $L^1(\Om)$ does not continuously embed into $(W_0^{1, 2}(\Om))^\star$
and that hence the Navier--Stokes energy identity does no longer provide an obvious starting point for bootstrap procedures.\abs
As observed in \cite[Theorem~1.2]{win_lp}, this problem can be overcome by additionally requiring a bound for $\int_0^T \io |f| \ln(|f|+1)$ as well.
Although this additional assumption is rather mild, it turns out to be too strong for the system considered in the second part of the present paper (see below),
the main problem being that the bounds derived in \cite[Theorem~1.2]{win_lp} could depend in an unfavorable way on the space-time $L\log L$ norm of $f$.\abs
In contrast to this, our first main result establishes uniform $L^p$ bounds which may also depend on the space-time $L^2$ norm of $f$, but only raised to an arbitrary small power.
That is, while we require a bound for $f$ in a stronger topology compared to \cite[Theorem~1.2]{win_lp}, we are able to control the influence of this bound in a convenient quantitative manner.
To be precise, in Section~\ref{sec:nasto} we shall see the following.
\begin{theo}\label{theo1}
  Let $\Om\subset\R^2$ be a bounded domain with smooth boundary, and let $T\in (0,\infty)$,
  $f\in C^0(\bom\times [0,T);\R^2)$, 
  $v\in \bigcup_{\beta\in (\frac{1}{2},1)} C^0([0,T);D(A_2^\beta)) \cap C^{2,1}(\bom\times (0,T);\R^2)$
  and $Q\in C^{1,0}(\Om\times (0,T))$ be such that \eqref{1.1} holds
  and that \eqref{0v} is solved in the classical sense. 
  Then for any $p \in [1, \infty)$ and each $\theta>0$ there exists $K=K(p,\theta,f,v)>0$ such that
  \be{1.01}
	\|v(\cdot,t)\|_{L^p(\Om)}
	\le \bigg\{ \int_0^t \io |f|^2 + K \bigg\}^\theta
	\qquad \mbox{for all } t\in (0,T).
  \ee
\end{theo}
Similarly as in \cite{win_lp}, our proof starts by splitting $v$ into the solution $v_1$ of the inhomogeneous Stokes equations
and a solution $v_2$ of the Navier--Stokes equations containing certain error terms but not explicitly depending on $f$.
Testing the latter with $A_2^{2\beta} v_2$ for $\beta \in (0, \frac12)$ close to $\frac12$ (see Lemma~\ref{lem4})
and making use of the energy identity (cf.\ Lemma~\ref{lem3}) yield $L^p$ estimates for $v$ essentially depending on a quantity exponential in $\int_{t_0}^T \io |f| \ln(|f|+1)$.
For $t_0$ sufficiently close to $T$, we can then conclude \eqref{1.01} for $t \in (t_0, T)$, while the regularity of $v$ entails \eqref{1.01} for $t \in [0, t_0]$, see Lemma~\ref{lem5}.
Of course, this approach necessitates that the constant $K$ in \eqref{1.01} may depend on $v$ and not just on the parameters and some norm of the data.\Abs
{\bf Application: A chemotaxis--Navier--Stokes system with local sensing.} \quad
Next, we consider
\be{0}
	\left\{ \begin{array}{rcll}
	n_t + u\cdot\nabla n
	&=& \Delta \big(n c^{-\al}\big),
	\qquad & x\in \Omega, \ t>0, \\[1mm]
	c_t + u\cdot\nabla c
	&=& \Delta c -nc,
	\qquad & x\in \Omega, \ t>0, \\[1mm]
	u_t + (u\cdot\na) u
	&=& \Delta u + \nabla P + n\nabla\Phi,
	\qquad \nabla\cdot u=0,
	\qquad & x\in \Omega, \ t>0, \\[1mm]
	& & \hspace*{-29mm}
	\frac{\partial n}{\partial\nu} = \frac{\partial c}{\partial\nu}=0, \quad u=0,
	\qquad & x\in\pO, \ t>0, \\[1mm]
	& & \hspace*{-29mm}
	n(x,0)=n_0(x), \quad c(x,0)=c_0(x), \quad u(x,0)=u_0(x),
	\qquad & x\in\Omega,
	\end{array} \right.
\ee
in a smoothly bounded domain $\Om\subset\R^2$, where $\al>0$ is arbitrary and $\Phi \in W^{2, \infty}(\Om)$ is a given function.
The system \eqref{0} models the behavior of bacteria (with density $n$) which may partially orient their movement towards oxygen (with density $c$)
and who interact with a fluid (with velocity $u$) through buoyancy and transportation.
The right-hand side of the first equation in \eqref{0} expands to $\nabla \cdot (D(c) \nabla n - n \chi(c) \nabla c)$ with $D(c) = c^{-\alpha}$ and $\chi(c) = \alpha c^{-\alpha-1}$, $c > 0$,
and for general choices of $D$ and $\chi$ such a system (without fluid) has already been proposed by Keller and Segel in the 1970s \cite{KellerSegelTravelingBandsChemotactic1971},
while the coupling to a Navier--Stokes system goes back to \cite{TuvalEtAlBacterialSwimmingOxygen2005}.\abs
Such systems have received considerable attention in the last decade;
especially for results regarding the fluid-free setting we refer to the recent survey \cite{LankeitWinklerDepletingSignalAnalysis2023}.
Here, we mainly confine ourselves with two-dimensional chemotaxis-fluid systems.
While in the arguably simplest case, namely when both $D$ and $\chi$ are positive constants,
global classical solutions always exist 	
(\cite{JiangEtAlGlobalExistenceAsymptotic2015}),
the situation becomes much more delicate when $\chi$ is singular at zero,
as the second solution component may become small and hence strengthen the potentially destabilizing chemotaxis term over time;
in fact, to underline this we recall that already in the classical Keller--Segel-production system
\bas
  \left\{ \begin{array}{rcl}
  n_t &=& \Delta n - \chi \nabla \cdot \big(n \nabla c \big), \\[1mm]
  c_t &=& \Delta c - c + n,
  \end{array} \right.
\eas
the size of $\chi$ (relative to $\io n_0$)
determines if all solutions remain bounded or if singularities can form for some initial data
(see, e.g., \cite{HorstmannWangBlowupChemotaxisModel2001}, \cite{NagaiEtAlApplicationTrudingerMoserInequality1997}, \cite{HerreroVelazquezBlowupMechanismChemotaxis1997}).\abs
%
This challenge is reflected by the fact that unconditional global existence results of classical solutions (or, alternatively, proofs of finite-time blow-up)
of chemotaxis-consumption models with fluid and logarithmic sensitivity (i.e., $D \equiv 1$ and $\chi(c) = \frac1c$ for $c > 0$),
hence accounting for stimulus perception in accordance with the Weber--Fechner law (\cite{RosenSteadystateDistributionBacteria1978}),
are apparently yet unavailable.
Indeed, up to now only some global generalized solutions have been constructed (\cite{WangGlobalLargedataGeneralized2016}, \cite{LiuLargetimeBehaviorTwodimensional2021})
which become smooth after some finite time if the mass of the bacteria is sufficiently small (\cite{BlackEventualSmoothnessGeneralized2018}, \cite{LiuLargetimeBehaviorTwodimensional2021})
or if the chemical is consumed sufficiently slowly~(\cite{FuestChemotaxisFluidSystems2023}).
As by-products, these latter works also obtain global classical solvability under certain smallness conditions.\abs
In settings where the motion of the organism is starvation-driven, also $D(c)$ may drastically increase as $c$ approaches zero,
see \cite{FuEtAlStripeFormationBacterial2012}, \cite{LiuEtAlSequentialEstablishmentStripe2011}, \cite{ChoKimStarvationDrivenDiffusion2013} for recent modelling considerations.
Mathematically, one might hope that when also the regularizing effect of the diffusion term is enhanced for small $c$,
the (provable) regularity of solutions improves compared to settings where only $\chi$ is singular at zero.
On a technical level, this is in particular the case when $D' = -\chi$, i.e., when $\nabla \cdot (D(c) \nabla n - n \chi(c) \nabla c) = \Delta(n D(c))$,
as then techniques based on duality arguments 		
become available.
This is the setting we consider here and which has also been proposed in the modelling works referenced above.\abs
This hope is indeed justified:
Global classical solutions for the fluid-free counterpart of \eqref{0} have very recently been constructed in \cite{win_log_GNI} for all $\al > 0$,
while global weak solutions are known to exist also in the higher-dimensional setting (\cite{taowin_JDE}).
For precedents with non-singular motility, where $c^{-\al}$ in the first equation in \eqref{0} is replaced by $(c+\eps)^{-\al}$ for positive $\eps$,
we refer to, e.g., \cite{LiWinklerRelaxationKellerSegelconsumptionSystem2023} 
and \cite{win_lp}.
Let us also briefly mention that chemotaxis systems with local sensing have also been studied in various related situations,
see for instance 
\cite{LiZhaoGlobalBoundednessLarge2021}, \cite{LaurencotLongTermSpatial2023}, 
\cite{wang_liangchen},
\cite{li_wang},
\cite{WinklerQuantitativeStrongParabolic2022}, \cite{WinklerGlobalGeneralizedSolvability2023}
for the non-singular but potentially degenerate setting 
and \cite{FujieSenbaGlobalExistenceInfinite2022}, 
\cite{FujieJiangComparisonMethodsKellerSegeltype2021}, 	
\cite{BurgerEtAlDelayedBlowupChemotaxis2021}, \cite{JinWangCriticalMassKellerSegel2020}, 
\cite{lv_wenbin},
\cite{DesvillettesEtAlLogarithmicChemotaxisModel2019} 
among others for the case of signal production. \abs
In contrast to \cite{win_log_GNI}, however, the coupling to the Navier--Stokes equations in \eqref{0} requires to
control the terms stemming from the transportation term $u \cdot \nabla n$ in order to gain any useful information from duality arguments.
This eventually comes down to estimating $u$ in some $L^p$ norm,
but as apart from mass conservation no a~priori information for the first solution component appears to be available at this point --  not even a space-time $L\log L$ bound --
the $L^p$ estimates derived in \cite{win_lp} are not applicable.
Hence, we make use of Theorem~\ref{theo1} instead.
If we choose $\theta > 0$ therein sufficiently small, the right-hand side in \eqref{1.01} results in only a small power of $\int_0^T \io n^2$ which can be absorbed into a dissipative term (cf.\ Lemma~\ref{lem6}).\abs
After this crucial first step, the additional presence of a fluid only requires minor modifications of the ideas from \cite{win_log_GNI},
so that we are able to bootstrap the bounds obtained in Lemma~\ref{lem6} to estimates so strong that they exclude the possibility of finite-time blow-up.\abs
Our main results concerning (\ref{0}) then reads as follows.
\begin{theo}\label{theo13}
  Let $\Om\subset\R^2$ be a bounded domain with smooth boundary, let $\al>0$ and $\Phi\in W^{2,\infty}(\Om)$, and suppose that
  \be{init}
	\left\{ \begin{array}{l}
	n_0\in W^{1,\infty}(\Om)
	\quad \mbox{is nonnegative, \qquad that} \\[1mm]
	c_0\in W^{1,\infty}(\Omega)
	\quad \mbox{is positive in $\bom$, \qquad and that} \\[1mm]
	u_0\in D(A_2^\iota) \mbox{ with some } \iota \in (\frac{1}{2},1).
	\end{array} \right.
  \ee
  Then there exist functions
  \be{reg}
	\lbal
	n \in C^0(\bom\times [0,\infty)) \cap C^{2,1}(\bom\times (0,\infty)), \\[1mm]
	c\in \bigcap_{q>2} C^0([0,\infty);W^{1,q}(\Om)) \cap C^{2,1}(\bom\times (0,\infty)), \\[1mm]
	u\in C^0([0,\infty);D(A_2^\iota)) \cap C^{2,1}(\bom\times (0,\infty);\R^2) 
	\qquad \mbox{and} \\[1mm]
	P \in C^{1,0}(\Om\times (0,\infty))
	\ear
  \ee
  such that $n\ge 0$ and $c>0$ in $\bom\times [0,\infty)$, and that (\ref{0}) is solved in the classical sense.
\end{theo}
\mysection{Estimates for solutions to \eqref{0v}. Proof of Theorem~\ref{theo1}}\label{sec:nasto}
Throughout this section, for $p\in (1,\infty)$ we denote the Stokes operator on $L_\sigma^p(\Om) = \{\varphi \in L^p(\Om; \R^2) \mid \nabla \cdot \varphi = 0 \text{ in } \mc D(\Om)\}$ by $A_p$
and the Helmholtz projection from $L^p(\Om; \R^2)$ to $L_\sigma^p(\Om)$ by $\mc P_p$ for $p \in (1, \infty)$.
Since these operators coincide on $C_0^\infty(\Omega; \R^2)$ and thus on the intersection of their domains, we may often write $A$ and $\mc P$ instead of $A_p$ and $\mc P_p$, respectively, without specifying $p$.
Moreover, we henceforth suppose that the hypotheses of Theorem~\ref{theo1} hold, and set $v_0 \defs v(\cdot, 0)$.\abs
We first note that letting
\be{v12}
v_1(\cdot,t)\defs \ure^{-tA} v_0 + \int_0^t \ure^{-(t-s)A} \proj f(\cdot,s) ds
\qquad \mbox{and} \qquad
v_2(\cdot,t)\defs v(\cdot,t)-v_1(\cdot,t),
\qquad t\in [0,T),
\ee
by relying on \eqref{1.1} in quite a straightforward manner we can obtain bounds for $v_1$, the solution of a Stokes system with force term $f$.
\begin{lem}\label{lem2}
  Under the assumptions of Theorem~\ref{theo1}, we have
  \bas
	\sup_{t\in (0,T)} \|v_1(\cdot,t)\|_{L^p(\Om)} < \infty
	\qquad \mbox{for all $p\in [1,\infty)$.}
  \eas
\end{lem}
\proof
  This follows by applying standard smoothing estimates for the Stokes semigroup (\cite[p.\ 201]{giga1986}) to (\ref{v12}) in the considered two-dimensional framework,
  see \cite[Lemma~2.5]{WangXiangGlobalExistenceBoundedness2015} 
  for details.
\qed
Next, we record an estimate resulting from a well-known combination of the Navier--Stokes energy identity combined with the Trudinger--Moser inequality,
which is implicitly contained in \cite[Lemma~2.7]{win_SIMA} or \cite[Lemma~2.5]{win_lp}, for instance.
Unlike in these works, however, we do not have a space-time $L\log L$ bound for $f$ at our disposal,
and hence the bounds in the following lemma do not yet provide any unconditional estimates but need to depend on $f$.
(That is, the right-hand side in \eqref{3.1} below may be unbounded for $t \nearrow T$.)
\begin{lem}\label{lem3}
  There exists $K_1>0$ such that for each $t_0\in [0,T)$ it is possible to fix $k_1(t_0)>0$ in such a way that
  \be{3.1}
	\int_{t_0}^t \io |\na v|^2 
	\le K_1 \int_{t_0}^t \io |f|\ln (|f|+1)
	+ k_1(t_0)
	\qquad \mbox{for all } t\in (t_0,T).
  \ee
\end{lem}
\proof
  In line with (\ref{1.1}), we fix $C_1>0$ such that
  \bas
	\io |f(\cdot,t)| \le C_1
	\qquad \mbox{for all } t\in (0,T),
  \eas
  and according to a Poincar\'e inequality we can find $C_2>0$ fulfilling
  \be{3.3}
	\io |\vp|^2 \le C_2 \io |\na\vp|^2
	\qquad \mbox{for all } \vp\in W_0^{1,2}(\Om;\R^2).
  \ee
  We next rely on a consequence of the Moser--Trudinger inequality in the planar domain $\Om$
  (see \cite[Lemma~2.2]{win_SIMA}) to pick $C_3>0$ in such a way that
  \bea{3.4}
	\io \rho_1\cdot\rho_2
	&\le& \frac{1}{4C_1} \cdot \bigg\{ \io |\rho_1|\bigg\} \io |\na\rho_2|^2
	+ \frac{1}{4C_1 C_2 |\Om|} \cdot \bigg\{ \io |\rho_1|\bigg\} \cdot \bigg\{ \io |\rho_2| \bigg\}^2 \nn\\
	& & + C_3 \io |\rho_1| \ln \big( |\rho_1|+1\big) 
	+ C_3 \cdot \bigg\{ \io |\rho_1| \bigg\} \cdot \left\{ 1 - \ln \bigg\{ \frac{1}{|\Om|} \io |\rho_1| \bigg\} \right\} \nn\\
	& & \mbox{whenever $0\not\equiv\rho_1 \in C^0(\bom;\R^2)$ and $\rho_2\in C^1(\bom,\R^2)$,}
  \eea
  so that with $C_4 \defs C_3 \big(C_1 + \frac{|\Omega|}{\ure}\big)$, thanks to the Cauchy--Schwarz inequality we obtain that
  \bas
	\io f\cdot v
	&\le& \frac{1}{4C_1} \cdot \bigg\{ \io |f| \bigg\} \cdot \io |\na v|^2
	+ \frac{1}{4C_1 C_2 |\Om|} \cdot \bigg\{ \io |f|\bigg\} \cdot \bigg\{ \io |v| \bigg\}^2 \nn\\
	& & + C_3 \io |f|\ln \big(|f|+1\big)
	+ C_3 \cdot \bigg\{ \io |f|\bigg\} \cdot \left\{ 1 - \ln \bigg\{ \frac{1}{|\Om|} \io |f|\bigg\} \right\}\nn\\
	&\le& \frac{1}{4} \io |\na v|^2
	+ \frac{1}{4C_2 |\Om|} \cdot \bigg\{ \io |v|\bigg\}^2
	+ C_3 \io |f|\ln \big(|f|+1\big) 
	+ C_4 \nn\\
	&\le& \frac{1}{4} \io |\na v|^2
	+ \frac{1}{4C_2} \io |v|^2 
	+ C_3 \io |f|\ln \big(|f|+1\big) 
	+ C_4 \nn\\
	&\le& \frac{1}{2} \io |\na v|^2
	+ C_3 \io |f|\ln \big(|f|+1\big) 
	+ C_4 
	\qquad \mbox{for all } t\in (0,T).
  \eas
  Therefore,
  \bas
	\frac{1}{2} \ddt \io |v|^2
	+ \io |\na v|^2
	= \io f\cdot v
	\le \frac{1}{2} \io |\na v|^2
	+ C_3 \io |f|\ln \big(|f|+1\big) 
	+ C_4 
	\qquad \mbox{for all } t\in (0,T)
  \eas
  and hence
  \bas
	\ddt \io |v|^2
	+ \io |\na v|^2
	\le 2C_3 \io |f|\ln \big(|f|+1\big) 
	+ 2 C_4 
	\qquad \mbox{for all } t\in (0,T),
  \eas
  which implies
  \bas
	\int_{t_0}^t \io |\na v|^2 \le 
	\io |v_1(\cdot,t_0)|^2 + 2C_4 t 
	+ 2C_3 \int_{t_0}^t \io |f|\ln \big(|f|+1\big)
  \eas
	for all $t_0\in [0,T)$ and $t\in (t_0,T)$
  and thus (\ref{3.1}) upon evident choices of $K_1$ and $(k_1(t_0))_{t_0\in (0,T)}$.
\qed
As used multiple times in the sequel, we briefly state embedding results regarding the domains of fractional powers of the Stokes operator,
including the case of zeroth power, that is, embeddings into $D(A_p^0) = L_\sigma^p(\Om)$.
\begin{lem}\label{lem3.5}
  Let $\beta_1, \beta_2 \in [0, 1]$ and $p_1, p_2 \in (1, \infty)$ with $2\beta_2 \neq \frac{1}{p_2}$, $\beta_1 \ge \beta_2$ and $2\beta_1 - \frac{2}{p_1} \ge 2\beta_2 - \frac{2}{p_2}$.
  Then $D(A_{p_1}^{\beta_1}) \hookrightarrow D(A_{p_2}^{\beta_2})$.
\end{lem}
\proof
  For $\beta \in [0, 1]$ and $p \in (1, \infty)$,
  we set $W_{\mc B}^{2\beta, p}(\Om) \defs \{\varphi \in W^{2\beta, p} \mid \varphi=0 \text{ on $\partial\Om$ if $2\beta > \frac1p$}\}$ 
  and introduce the operator $\Lambda_p \defs -\Delta$ on $L^p(\Om)$ with $D(\Lambda_p) = W_{\mc B}^{2, p}(\Om)$.
  According to \cite[Theorems~1.15.3 and 4.3.3]{triebel}, we then have $D(\Lambda_p^\beta) = [L^p(\Om), W_{\mc B}^{2, p}(\Om)]_\beta \subseteq W_{\mc B}^{2\beta,p}(\Om)$ for $\beta \in [0, 1]$ and $p \in (1, \infty)$,
  where equality holds if $2\beta \neq \frac{1}{p}$.
  Since $W^{2\beta_1, p_1}(\Om) \hookrightarrow W^{2\beta_2, p_2}(\Omega)$ by \cite[Theorem~4.6.2]{triebel},
  we conclude $D(\Lambda_{p_1}^{\beta_1}) \hookrightarrow D(\Lambda_{p_2}^{\beta_2})$,
  upon which the statement results by \cite[Theorem~3]{giga1981_the_other}. 
\qed
With Lemma~\ref{lem3.5} at hand, we may now test the equation solved by $v_2$ -- which, importantly, does not explicitly depend on the force $f$ --
with $A^{2\beta} v_2$ to obtain uniform-in-time $L^p$ estimates also for $v_2$.
However, they yet depend on the quantity in the left-hand side of \eqref{3.1} and will then be combined with Lemma~\ref{lem3} in Lemma~\ref{lem5} below.
As we shall see in the proof of the latter, it turns out to be crucial that the constant $K_2$ appearing in \eqref{4.1} below does not depend on $t_0$.
\begin{lem}\label{lem4}
  For each $p\ge 1$ there exists $K_2(p)>0$ with the property that whenever 
  $t_0\in [0,T)$, one can fix $k_2(p,t_0)>0$ in such a way that
  \be{4.1}
	\|v(\cdot,t)\|_{L^p(\Om)} \le k_2(p,t_0) \ure^{K_2(p) \int_{t_0}^t \io |\na v|^2}
	\qquad \mbox{for all } t\in (t_0,T).
  \ee
\end{lem}
\proof
  Lemma~\ref{lem3.5} allows us to choose $\beta=\beta(p)\in (0,\frac{1}{2})$ suitably close to $\frac{1}{2}$ such that with some $C_1=C_1(p)>0$ we have
  \be{4.2}
	\|\vp\|_{L^p(\Om)} \le C_1 \|A^\beta \vp\|_{L^2(\Om)}
	\qquad \mbox{for all } \vp\in D(A^\beta),
  \ee
  and then obtain on testing the identity $\pa_t v_2+(v\cdot\na) v=\Del v_2 + \na Q_2$, satisfied with some
  $Q_2\in C^{1,0}(\Om\times (0,T))$ according to (\ref{0v}) and (\ref{v12}), by $A^{2\beta} v_2$ that
  \be{4.3}
	\ddt \io |A^\beta v_2|^2
	+ 2 \io \big| A^\frac{2\beta+1}{2} v_2\big|^2
	= - 2 \io A^{2\beta} v_2 \cdot \proj \big[ (v\cdot\na)v\big]
	\qquad \mbox{for all } t\in (0,T).
  \ee
  Here, since $D(A_2^\frac{2\beta+1}{2}) \hra D(A_\frac{1}{\beta}^{2\beta})$ by Lemma~\ref{lem3.5}
  and as $\proj$ is continuous on $L^\frac{1}{1-\beta}(\Om;\R^2)$ (\cite[Theorem~1]{fujiwara_morimoto}), there exist $C_2=C_2(p)>0$ and
  $C_3=C_3(p)>0$ such that due to the H\"older inequality we have
  \bas
	- 2 \io A^{2\beta} v_2 \cdot \proj \big[ (v\cdot\na)v\big]
	&\le& \|A^{2\beta} v_2\|_{L^\frac{1}{\beta}(\Om)} \big\| \proj[(v\cdot\na) v]\big\|_{L^\frac{1}{1-\beta}(\Om)} \\
	&\le& C_2 \|A^\frac{2\beta+1}{2} v_2\|_{L^2(\Om)} \|(v\cdot\na) v\|_{L^\frac{1}{1-\beta}(\Om)} \\
	&\le& 2 \io \big| A^\frac{2\beta+1}{2} v_2\big|^2 \\
	& & + C_3 \|(v_1\cdot\na)v\|_{L^\frac{1}{1-\beta}(\Om)}^2
	+ C_3 \|(v_2\cdot\na)v\|_{L^\frac{1}{1-\beta}(\Om)}^2
	\qquad \mbox{for all } t\in (0,T),
  \eas
  where again by the H\"older inequality,
  \bas
	\|(v_i\cdot\na)v\|_{L^\frac{1}{1-\beta}(\Om)}^2
	\le \|v_i\|_{L^\frac{1-2\beta}{1-\beta}(\Om)}^2 \|\na v\|_{L^2(\Om)}^2
	\qquad \mbox{for $i\in \{1,2\}$ and all } t\in (0,T).
  \eas
  Since $\beta < \frac12$, employing Lemma~\ref{lem2} we thus find $C_4=C_4(p)>0$ such that
  \bas
	\|(v_1\cdot\na)v\|_{L^\frac{1}{1-\beta}(\Om)}^2
	\le C_4 \io |\na v|^2
	\qquad \mbox{for all } t\in (0,T),
  \eas
  while by Hölder's inequality and since $D(A_2^\beta) \hra L^\frac{2}{1-2\beta}(\Om)$ according to Lemma~\ref{lem3.5} (which is applicable as $\beta < \frac12$),
  we can fix $C_5=C_5(p)>0$ satisfying
  \bas
	\|(v_2\cdot\na)v\|_{L^\frac{1}{1-\beta}(\Om)}^2
	\le C_5 \cdot \bigg\{ \io |\na v|^2 \bigg\} \cdot \io |A^\beta v_2|^2
	\qquad \mbox{for all } t\in (0,T).
  \eas
  Therefore, (\ref{4.3}) entails that $y(t)\defs \io |A^\beta v_2(\cdot,t)|^2 + 1$, $t\in [0,T)$, 
  and $C_6= C_6(p)\defs C_3(C_4+C_5)$
  have the property that
  \bas
	y'(t) \le C_6 \cdot \bigg\{ \io |\na v|^2 \bigg\} \cdot y(t)
	\qquad \mbox{for all } t\in (0,T),
  \eas
  and that thus
  \bas
	y(t) \le \bigg\{ \io |A^\beta v_2(\cdot,t_0)|^2 + 1 \bigg\} \cdot \ure^{C_6 \int_{t_0}^t \io |\na v|^2}
	\qquad \mbox{for all $t_0\in [0,T)$ and each } t\in (t_0,T),
  \eas
  which in view of (\ref{4.2}) and Lemma~\ref{lem2} establishes the claim.
\qed
Estimates of the form $\|v(\cdot, t)\|_{L^p(\Om)} \le k(t_0) \ure^{K \int_{t_0}^t \io |f| \ln(|f|+1)}$ are direct consequences of Lemma~\ref{lem3} and Lemma~\ref{lem4}.
At least on small time scales, i.e., if $t_0$ is sufficiently close to $T$, we can estimate the right-hand side therein against a small power of the space-time $L^2$ norm of $f$.
Since $K$ does not depend on $t_0$, we can conclude
\begin{lem}\label{lem5}
  For any $p\ge 1$ and each $\theta>0$ there exists $K>0$ such that
  \bas
	\|v(\cdot,t)\|_{L^p(\Om)}
	\le \bigg\{ \int_0^t \io |f|^2 + K \bigg\}^\theta
	\qquad \mbox{for all } t\in (0,T).
  \eas
\end{lem}
\proof
  According to \cite[Lemma 3.6]{win_NoDEA}, we have
  \bas
	\io \vp\ln\vp \le \bigg\{ \io \vp + 1 \bigg\} \cdot \ln \bigg\{ \io \vp^2 + \ure \bigg\}
	\qquad \mbox{for all nonnegative } \vp\in C^0(\bom),
  \eas
  so that by Young's inequality,
  \bea{1.2}
	\io |f|\ln (|f|+1)
	&\le& \io (|f|+1) \ln (|f|+1) \nn\\
	&\le& \bigg\{ \io (|f|+1) + 1 \bigg\} \cdot \ln \bigg\{ \io (|f|+1)^2 + \ure\bigg\} \nn\\
	&\le& C_1 \ln \bigg\{ 2\io |f|^2 + C_2 \bigg\}
	\qquad \mbox{for all } t\in (0,T),
  \eea
  where $C_2\defs 2|\Om|+\ure$, and where $C_1\defs \sup_{t\in (0,T)} \io \big( |f(\cdot,t)|+1\big) +1$ is finite thanks to (\ref{1.1}).\\
  We now fix $t_0\in [0,T)$ sufficiently close to $T$ such that with $K_1$ and $K_2$ taken from Lemma~\ref{lem3} and Lemma~\ref{lem4}
  we have
  \be{1.3}
	K_1 K_2 C_1\cdot (T-t_0) \le \theta,
  \ee
  and combine the latter two lemmata with (\ref{1.2}) to see that with $k_1=k_1(t_0)$ and $k_2=k_2(p,t_0)$ introduced there we have
  \bas
	\|v(\cdot,t)\|_{L^p(\Om)}
	&\le& k_2 \ure^{K_2 \int_{t_0}^t \io |\na v|^2} \\
	&\le& k_2 \ure^{K_1 K_2 \int_{t_0}^t \io |f|\ln (|f|+1) + k_1 K_2} \\
	&\le& k_2 \ure^{k_1 K_2} \ure^{K_1 K_2 C_1 \int_{t_0}^t \ln \big\{ 2\io |f|^2 + c_2\big\}}
	\qquad \mbox{for all } t\in (t_0,T).
  \eas
  As the Jensen inequality together with the fact that $\xi\ln\frac{1}{\xi} \le \frac{1}{\ure}$ for all $\xi>0$ implies that here
  \bas
	K_1 K_2 C_1 \int_{t_0}^t \ln \bigg\{ 2\io |f|^2 + C_2 \bigg\} 
	&\le& K_1 K_2 C_1\cdot (t-t_0) \cdot \ln \Bigg\{ \frac{1}{t-t_0} \int_{t_0}^t \bigg\{ 2\io |f|^2+C_2 \bigg\} \Bigg\} \\
	&=& K_1 K_2 C_1\cdot (t-t_0) \cdot \ln \Bigg\{ \int_{t_0}^t \io |f|^2 + \frac{C_2}{2} \cdot (t-t_0)\Bigg\} \\
	& & + K_1 K_2 C_1\cdot (t-t_0) \cdot\ln \frac{2}{t-t_0} \\
	&\le& K_1 K_2 C_1\cdot (t-t_0) \cdot \ln \Bigg\{ \int_{t_0}^t \io |f|^2 + \frac{C_2 T}{2} +1 \Bigg\}
	+ \frac{2 K_1 K_2 C_1}{\ure} \\
	&\le& \theta \ln \bigg\{ \int_{t_0}^t \io |f|^2 + \frac{C_2 T}{2} + 1 \bigg\}
	+ \frac{2 K_1 K_2 C_1}{\ure} 
  \eas
	for all $t\in (t_0,T)$ due to (\ref{1.3}), we therefore obtain that
  \bas
	\|v(\cdot,t)\|_{L^p(\Om)}
	&\le& k_2 \ure^{k_1 K_2} \ure^{\frac{2K_1 K_2 C_1}{\ure}} \cdot \bigg\{ \int_0^t \io |f|^2 + \frac{C_2 T}{2} +1 \bigg\}^\theta
	\qquad \mbox{for all } t\in (t_0,T).
  \eas
  Since $v$ is continuous and hence bounded in $\bom\times [0,t_0]$, this yields (\ref{1.01}) with some suitably large
  $K=K(p,\theta,f,v)>0$.
\qed
\proofc of Theorem~\ref{theo1}. \quad
  The desired estimate has been proven in Lemma~\ref{lem5}.
\qed
\mysection{Analysis of (\ref{0}). Proof of Theorem~\ref{theo13}}
We begin our study of \eqref{0} by stating a local existence result as well as the basic bounds for $n$ and $c$ in $L^\infty$-$L^1$ and $L^\infty$-$L^\infty$, respectively.
\begin{lem}\label{lem_loc}
  Let $\Om\subset\R^2$ be a bounded domain with smooth boundary, let $\al>0$ and $\Phi\in W^{2,\infty}(\Om)$, and assume
  (\ref{init}).
  Then there exist $\tm\in (0,\infty]$ as well as functions
  \bas
	\lbal
	n \in C^0(\bom\times [0,\tm)) \cap C^{2,1}(\bom\times (0,\tm)), \\[1mm]
	c\in \bigcap_{q>2} C^0([0,\tm);W^{1,q}(\Om)) \cap C^{2,1}(\bom\times (0,\tm)), \\[1mm]
	u\in C^0([0,\tm);D(A^\iota)) \cap C^{2,1}(\bom\times (0,\tm);\R^2) 
	\qquad \mbox{and} \\[1mm]
	P \in C^{1,0}(\Om\times (0,\tm))
	\ear
  \eas
  such that $n\ge 0$ and $c>0$ in $\bom\times [0,\tm)$, that $(n,c,u,P)$ solves (\ref{0}) classically in 
  $\Om\times (0,\tm)$, and that
  \be{ext}
	\mbox{if $\tm<\infty$, \quad then \quad}
	\limsup_{t\nearrow\tm} \Big\{ \|n(\cdot,t)\|_{L^\infty(\Om)}
	+ \|c(\cdot,t)\|_{W^{1,\infty}(\Om)}
	+ \|A^\iota u(\cdot,t)\|_{L^2(\Om)} \Big\}
	=\infty.
  \ee
  Moreover,
  \be{mass}
	\io n(\cdot,t)=\io n_0
	\qquad \mbox{for all } t\in (0,\tm)
  \ee
  as well as
  \be{cinfty}
	\|c(\cdot,t)\|_{L^\infty(\Om)} \le \|c_0\|_{L^\infty(\Om)}
	\qquad \mbox{for all } t\in (0,\tm).
  \ee
\end{lem}
\proof
  This can be seen by adapting a standard contraction mapping argument (\cite[Lemma~2.1]{JinDCDS2018}) to the present situation.
  (We note that as long $n$ is bounded, $c$ is at least locally in time bounded from below by a positive constant, so that we do not need to include a term such as $\|\frac{1}{c(t)}\|_{L^\infty(\Om)}$ in \eqref{ext}.)
\qed
Without commenting on this explicitly any further, we shall below assume that the smoothly bounded domain $\Om\subset\R^2$,
the number $\al>0$ and the function $\Phi\in W^{2,\infty}(\Om)$ have been fixed,
and let $(n,c,u,P)$ and $\tm$ be as thereupon obtained in Lemma~\ref{lem_loc}.\abs
The most crucial part of our analysis is performed in the following lemma.
Similarly as for instance in \cite{taowin_JDE}, \cite{tian_win} or \cite{win_log_GNI},
the basic idea is to make use of the structure of the first equation in \eqref{0} which allows for reasonings based on duality arguments.
The key additional challenge lies in the fact that no helpful a~prioi estimates for the fluid equation appear to be available only based on \eqref{mass}.
(This stays in contrast to situations where at least space-time $L\log L$ bounds for $n$ can be obtained in a rather straightforward manner, such as in \cite[Lemma~3.2]{win_lp}.)
Instead, we rely on Theorem~\ref{theo1} and estimate terms stemming from the fluid transportation term $u \cdot \nabla n$ essentially against the space-time $L^2$ norm of $n$,
whose boundedness is not established prior to this lemma but which appears as a dissipative term of the functional considered.
\begin{lem}\label{lem6}
  If $\tm<\infty$, then there exists $C>0$ such that
  \be{6.1}
	\int_0^t \io n^2 c^{-\al} \le C
	\qquad \mbox{for all } t\in (0,\tm)
  \ee
  and
  \be{6.01}
	\io c^{-\al}(\cdot,t) \le C
	\qquad \mbox{for all } t\in (0,\tm).
  \ee
\end{lem}
\proof
  Letting $B$ denote the realization of $-\Del$ under homogeneous Neumann boundary conditions in 
  $L^2_\perp(\Om)\defs \{\vp\in L^2(\Om) \mid \io \vp=0\}$, we follow \cite[Lemma 3.1]{taowin_JDE} and use the first two equations
  in (\ref{0}) to see that
  \be{6.2}
	\frac{1}{2} \ddt \io |B^{-\frac{1}{2}}(n-\onz)|^2
	+ \io n^2 c^{-\al}
	= \onz \io n c^{-\al}
	- \io (n-\onz) B^{-1} \na \cdot (nu)
	\qquad \mbox{for all } t\in (0,\tm),
  \ee
  and that
  \be{6.3}
	\frac{1}{\al} \ddt \io c^{-\al} + (\al+1) \io c^{-\al-2} |\na c|^2 = \io n c^{-\al}
	\qquad \mbox{for all } t\in (0,\tm),
  \ee
  where by Young's inequality,
  \be{6.4}
	(\onz+1) \io n c^{-\al} 
	\le \frac{1}{4} \io n^2 c^{-\al}
	+ (\onz+1)^2 \io c^{-\al}
	\qquad \mbox{for all } t\in (0,\tm).
  \ee
  Abbreviating $C_1 \defs  (\onz+1)^2\alpha > 0$, we thus see that the function $y \defs  \frac{1}{2} \io \big|B^{-\frac{1}{2}}(n-\onz)\big|^2 + \frac{1}{\al} \io c^{-\al}$ fulfills
  \bas
	y'(t) \le C_1 y(t) - \frac{3}{4} \io n^2 c^{-\al} - \io (n-\onz) B^{-1} \na \cdot (nu)
	\qquad \mbox{for all } t\in (0,\tm).
  \eas
  By the variations-of-constants formula, this implies
  \bea{eq:y}
  y(t)
  &\le& \ure^{C_1 t} y(0) - \frac34 \int_0^t \io \ure^{C_1(t-s)} n^2 c^{-\al} - \int_0^t \io \ure^{C_1(t-s)} (n-\onz) B^{-1} \na \cdot (nu) \notag \\
  &\le& C_2 y(0) - \frac34 \int_0^t \io n^2 c^{-\al} + C_2 \int_0^t \io \big|(n-\onz) B^{-1} \na \cdot (nu)\big|
  \eea
  for all $t \in (0, \tm)$, where $C_2 \defs  \ure^{C_1 \tm} < \infty$ and where $y(0)$ is finite by \eqref{init}.
  To estimate the rightmost summand in (\ref{eq:y}), following \cite[Lemma 4.3]{tian_win} we fix any $p\in (1,\frac{4}{3})$ and
  $\theta\in (0,\frac{2p-2}{p})$ and use that then an argument based on elliptic regularity (\cite[Lemma 4.1]{tian_win})
  provides $C_3>0$ such that
  \bas
	\|B^{-1} \na\cdot\vp\|_{L^\frac{2p}{2-p}(\Om)} 
	\le C_3 \|\vp\|_{L^p(\Om)}
	\qquad \mbox{for all $\vp\in C^1(\bom;\R^2)$ such that $\vp\cdot\nu=0$ on } \pO.
  \eas
  Therefore, using the H\"older inequality, Theorem~\ref{theo1} and (\ref{mass}) we obtain $C_4>0$ and $C_5>0$ such that
  \bas
	\io \big|(n-\onz) B^{-1} \na\cdot (nu) \big|
	&\le& \|n-\onz\|_{L^\frac{2p}{3p-2}(\Om)} \|B^{-1} \na\cdot (nu)\|_{L^\frac{2p}{2-p}(\Om)} \\
	&\le& C_3 \|n-\onz\|_{L^\frac{2p}{3p-2}(\Om)} \|nu\|_{L^p(\Om)} \\
	&\le& C_3 \|n-\onz\|_{L^\frac{2p}{3p-2}(\Om)} \|n\|_{L^\frac{2p}{3p-2}(\Om)} \|u\|_{L^\frac{2p}{4-3p}(\Om)} \\
	&\le& C_4 \|n\|_{L^\frac{2p}{3p-2}(\Om)}^2 \cdot \bigg\{ \int_0^t \io n^2 + C_4 \bigg\}^\theta \\
	&\le& C_4 \|n\|_{L^2(\Om)}^\frac{2(2-p)}{p} \|n\|_{L^1(\Om)}^\frac{4(p-1)}{p} 
		\cdot \bigg\{ \int_0^t \io n^2 + C_4 \bigg\}^\theta \\
	&\le& C_5 \|n\|_{L^2(\Om)}^\frac{2(2-p)}{p} \cdot \bigg\{ \int_0^t \io n^2 + C_4 \bigg\}^\theta
	\qquad \mbox{for all } t\in (0,\tm),
  \eas
  so that since
  \bas
	\io n^2 \le \|c_0\|_{L^\infty(\Om)}^\al \io n^2 c^{-\al}
	\qquad \mbox{for all } t\in (0,\tm)
  \eas
  due to (\ref{cinfty}), by Young's inequality it follows that with some $C_6>0$ and $C_7>0$,
  \bas
	\io \big|(n-\onz) B^{-1} \na\cdot (nu)\big|
	&\le& C_6 \cdot \bigg\{ \io n^2 c^{-\al} \bigg\}^\frac{2-p}{p} \cdot \bigg\{ \int_0^t \io n^2 c^{-\al} + C_6 \bigg\}^\theta \\
	&\le& \frac{1}{4C_2} \io n^2 c^{-\al} + C_7 \cdot \bigg\{ \int_0^t \io n^2 c^{-\al} \bigg\}^\frac{p\theta}{2p-2} + C_7
  \eas
	for all $t\in (0,\tm)$.
  Again making use of Young's inequality, we further obtain
  \bas
	\int_0^t \bigg\{ \int_0^s \io n^2(x, \sigma) c^{-\al}(x, \sigma) \dx \dsigma \bigg\}^\frac{p\theta}{2p-2} \ds
  &\le& \int_0^t \bigg\{ \int_0^t \io n^2(x, \sigma) c^{-\al}(x, \sigma) \dx \dsigma \bigg\}^\frac{p\theta}{2p-2} \ds \\
  &=& (4C_2C_7)^\frac{p\theta}{2p-2} t \cdot \bigg\{ \frac{1}{4C_2C_7} \int_0^t \io n^2 c^{-\al} \bigg\}^\frac{p\theta}{2p-2} \\
  &\le& \frac{1}{4C_2C_7} \int_0^t \io n^2 c^{-\al} + \Big\{ (4C_2C_7)^\frac{p\theta}{2p-2} \tm \Big\}^\frac{2p-2}{2p-2-p\theta}
  \eas
	for all $t\in (0,\tm)$, because $\frac{p\theta}{2p-2} < 1$.
  This first implies
  \bas
  \int_0^t \io \big|(n-\onz) B^{-1} \na\cdot (nu)\big| 
	\le \frac{1}{2C_2} \int_0^t \io n^2 c^{-\al} + C_8
  \eas
	for all $t\in (0,\tm)$ and some $C_8 > 0$, which is finite due to our assumption that $T_{\max}$ is finite.
  Combined with \eqref{eq:y}, we conclude
  \bas
  y(t) + \frac14\int_0^t \io n^2 c^{-\al}
  \le C_2 y(0) + C_2 C_8
	\qquad \mbox{for all } t\in (0,\tm),
  \eas
  so that by the definition of $y$ this establishes both (\ref{6.1}) and (\ref{6.01}) with some $C>0$.
\qed
With \eqref{6.1} and \eqref{6.01} at hand, the global existence proof proceeds quite analogously to the fluid-free setting considered in \cite[Lemma 3.3]{win_log_GNI}.
However, for the sake of completeness and as some modifications are necessary, we choose to at least sketch the proofs for the bootstrap procedure.
\begin{lem}\label{lem7}
  Assume that $\tm<\infty$. Then one can find $C>0$ such that
  \be{7.1}
	\int_0^t \io \frac{|\na c|^4}{c^3} \le C
	\qquad \mbox{for all } t\in (0,\tm)
  \ee
  and
  \be{7.2}
	\int_0^t \io \frac{n}{c} |\na c|^2 \le C
	\qquad \mbox{for all } t\in (0,\tm).
  \ee
\end{lem}
\proof
  A standard argument using the second equation in (\ref{0}) and the solenoidality of $u$ 
  (\cite[Lemma 3.3]{win_log_GNI})
  provides $C_1>0$ and $C_2>0$ such that
  \bas
	\ddt \io \frac{|\na c|^2}{c} 
	+ C_1 \io \frac{|\na c|^4}{c^3}
	+ C_1 \io \frac{n}{c} |\na c|^2
	\le C_2 \io n^2
	+ C_2 \io |u|^4
	\qquad \mbox{for all } t\in (0,\tm).
  \eas
  Therefore, 
  \bas
	C_1 \int_0^t \io \frac{|\na c|^4}{c^3} 
	+ C_1 \int_0^t \io \frac{n}{c} |\na c|^2
	\le \io \frac{|\na c_0|^2}{c_0}
	+ C_2 \int_0^t \io n^2
	+ C_2 \int_0^t \io |u|^4
  \eas
	for all $t\in (0,\tm)$,
  so that since Theorem~\ref{theo1} yields $C_3>0$ fulfilling
  \bas
	C_2 \int_0^t \io |u|^4
	&\le& C_2 \tm \cdot \sup_{t\in (0,\tm)} \|u(\cdot,t)\|_{L^4(\Om)}^4 \\
	&\le& C_3 \cdot \bigg\{ \int_0^{\tm} \io n^2 +1 \bigg\}
	\qquad \mbox{for all } t\in (0,\tm),
  \eas
  and since $\int_0^{\tm} \io n^2$ is finite by Lemma~\ref{lem6} and (\ref{cinfty}), from this we already obtain
  (\ref{7.1}) and (\ref{7.2}).
\qed
\begin{lem}\label{lem8}
  Let $\ell\in C^2((0,\infty))$ be such that $\xi\ell''(\xi)+2\ell'(\xi)\ge 0$ for all $\xi>0$. Then $z\defs nc^{-\al}$ satisfies
  \bea{8.1}
	& & \hs{-20mm}
	\ddt \io n\ell(z)
	+ \frac{1}{2} \io \big\{ z\ell''(z) + 2\ell'(z)\big\} \cdot |\na z|^2 \nn\\
	&\le& \al(\al-1) \io nc^{-2} \cdot z\ell'(z) |\na c|^2 \nn\\
	& & + \frac{\al^2}{2} \io c^{2\al-2} \cdot \big\{ z^3\ell''(z) + 2z^2 \ell'(z)\big\} \cdot |\na c|^2 \nn\\
	& & + \al \io n^2 \cdot z\ell'(z)
	\qquad \mbox{for all } t\in (0,\tm).
  \eea
\end{lem}
\proof
  Since the divergence-free vector field $u$ vanishing on $\pO$ does effectively not influence the evolution of
  $[0,\tm) \ni t\mapsto \io n(\cdot,t)\ell(z(\cdot,t))$, this can be derived by a verbatim copy of \cite[Lemma 3.5]{win_log_GNI}.
\qed
\begin{lem}\label{lem9}
  If $\tm<\infty$, then there exists $C>0$ such that
  \be{9.1}	
	\int_0^t \io n^2 \ln (n+\ure) \le C
	\qquad \mbox{for all } t\in (0,\tm).
  \ee
\end{lem}
\proof
  As detailed in \cite[Lemma 3.6 and Lemma 3.7]{win_log_GNI}, this can be seen on choosing $\ell(\xi)\defs \ln\xi, \xi>0$, in 
  Lemma~\ref{lem8}, and using that the bounds of the form
  \bas
	\sup_{t\in (0,\tm)} \io n(\cdot,t) \ln \big\{ n(\cdot,t)+\ure\big\} <\infty
  \eas
  and
  \bas
	\int_0^{\tm} \io \big|\na (nc^{-\al})\big|^2
	+ \int_0^{\tm} \io nc^{-2} |\na c|^2 < \infty,
  \eas
  as thereby obtained, imply (\ref{9.1}) through a Gagliardo--Nirenberg type inequality stated in
  \cite[Proposition 1.1]{win_log_GNI}.
\qed
\begin{lem}\label{lem10}
  Assume that $\tm<\infty$. Then there exist $\gamma>1$ and $C>0$ such that
  \be{10.1}	
	\int_0^t \io n^2 \ln^\gamma (n+\ure) \le C
	\qquad \mbox{for all } t\in (0,\tm).
  \ee
\end{lem}
\proof
  Relying on (\ref{9.1}) in estimating corresponding ill-signed terms on the right-hand side of (\ref{8.1}), this follows upon
  applying Lemma~\ref{lem8} to $\ell(\xi)\defs \ln^\gamma(\xi+\ure), \xi>0$, with arbitrary $\gamma\in (1,\frac{3}{2})$ if $\al\le 1$,
  and any $\gamma\in (1,\frac{3\al-1}{2\al})$ when $\al>1$, and again employing the interpolation result recorded in
  \cite[Proposition 1.1]{win_log_GNI} (cf.~\cite[Lemma~3.9 and Lemma~3.10]{win_log_GNI} for details).
\qed
\begin{lem}\label{lem11}
  If $\tm<\infty$, then there exists $C>0$ such that
  \be{11.1}	
	c(x,t)\ge C
	\qquad \mbox{for all $x\in\Om$ and } t\in (0,\tm).
  \ee
\end{lem}
\proof
  Modifying the argument from \cite[Lemma 3.11]{win_log_GNI}, we let $w\defs \ln \frac{\|c_0\|_{L^\infty(\Om)}}{c}$ and obtain that
  by (\ref{0}) and Young's inequality,
  \bas
	w_t &=& \Del w - |\na w|^2 + n - u\cdot\na w \\
	&\le& \Del w + n + \frac{1}{4} |u|^2
	\qquad \mbox{in } \Om\times (0,\tm).
  \eas
  Since $\frac{\pa w}{\pa\nu}=0$ on $\pO\times (0,\tm)$, according to a comparison principle this implies that
  \be{11.2}
	w\le \ow 
	\qquad \mbox{in } \Om\times (0,\tm),
  \ee
  where $\ow\in C^0(\bom\times [0,\tm)) \cap C^{2,1}(\bom\times (0,\tm))$ denotes the classical solution of the inhomogeneous
  linear problem
  \bas
	\lball
	\ow_t = \Del\ow + n + \frac{1}{4} |u|^2,
	\qquad & x\in \Omega, \ t\in (0,\tm), \\[1mm]
	\frac{\pa \ow}{\pa\nu}=0,
	\qquad & x\in \pO, \ t\in (0,\tm), \\[1mm]
	\ow(x,0)=\ln\frac{\|c_0\|_{L^\infty(\Om)}}{c_0(x)},
	\qquad & x\in\Om.
	\ear
  \eas
  But from Lemma~\ref{lem10}, Theorem~\ref{theo1} and Lemma~\ref{lem6} it readily follows that with some $\gamma>1$
  and $C_1>0$ we have
  \bas
	\int_0^t \io \Big\{ n + \frac{1}{4} |u|^2 \Big\}^2 \cdot \ln^\gamma \Big\{ n+ \frac{1}{4} |u|^2 + \ure \Big\} \le C_1
	\qquad \mbox{for all } t\in (0,\tm),
  \eas
  so that a recent result on regularity in heat equations involving sources from Orlicz spaces (\cite[Corollary~1.3]{win_IMRN})
  applies to ensure that since here we have $\gamma>1$, it follows that with some $C_2>0$,
  \bas
	\|\ow(\cdot,t)\|_{L^\infty(\Om)} \le C_2
	\qquad \mbox{for all } t\in (0,\tm).
  \eas
  Together with (\ref{11.2}) and the nonnegativity of $w$ asserted by (\ref{cinfty}), this implies (\ref{11.1}).
\qed
\begin{lem}\label{lem12}
  Suppose that $\tm<\infty$, and let $\iota$ be as in (\ref{init}). Then there exists $C>0$ such that
  \be{12.1}	
	\|n(\cdot,t)\|_{L^\infty(\Om)} \le C
	\qquad \mbox{for all } t\in (0,\tm),
  \ee
  that
  \be{12.2}	
	\|c(\cdot,t)\|_{W^{1,\infty}(\Om)} \le C
	\qquad \mbox{for all } t\in (0,\tm),
  \ee
  and that
  \be{12.3}	
	\|A^\iota u(\cdot,t)\|_{L^2(\Om)} \le C
	\qquad \mbox{for all } t\in (0,\tm).
  \ee 
\end{lem}
\proof
  For arbitrary $p>1$, straightforward manipulations applied to the first equation in (\ref{0}) (\cite[Lemma 3.12]{win_log_GNI})
  provide $C_1, C_2, C_3, C_4 >0$ such that in line with Lemma~\ref{lem11} we have
  \bas
	& & \hs{-20mm}
	\frac{1}{p} \ddt \io n^p
	+ C_1 \io \big|\na n^\frac{p}{2}\big|^2 \\
	&=& - (p-1) \io n^{p-2} c^{-\al} |\na n|^2
	+ (p-1)\al \io n^{p-1} c^{-\al-1} \na n\cdot\na c
	+ C_1 \io \big|\na n^\frac{p}{2}\big|^2 \\
	&\le& - \frac{C_1}{2} \io \big|\na n^\frac{p}{2}\big|^2
	+ \frac{(p-1)\al^2}{2} \io n^p c^{-\al-2} |\na c|^2 \\
	&\le& - \frac{C_1}{2} \io \big|\na n^\frac{p}{2}\big|^2
	+ C_2 \io n^p \big|\na c^\frac{1}{4}\big|^2 \\
	&\le& - \frac{C_1}{2} \io \big|\na n^\frac{p}{2}\big|^2
	+ C_3 \big\|\na c^\frac{1}{4}\big\|_{L^4(\Om)}^2 \big\|\na n^\frac{p}{2}\big\|_{L^2(\Om)} 
		\big\| n^\frac{p}{2}\big\|_{L^2(\Om)} \\
	& & + C_3\big\|\na c^\frac{1}{4}\big\|_{L^4(\Om)}^2 \big\|n^\frac{p}{2}\big\|_{L^2(\Om)}^2 \\
	&\le& C_4\cdot\bigg\{ \io \frac{|\na c|^4}{c^3} + 1 \bigg\} \cdot \io n^p
	\qquad \mbox{for all } t\in (0,\tm).
  \eas
  According to Lemma~\ref{lem7}, this implies that $\sup_{t\in(0,\tm)} \io n^p(\cdot,t) <\infty$ for all $p\in (1,\infty)$,
  which due to standard smoothing estimates for the Stokes and the heat semigroup firstly implies (\ref{12.3})
  and thereupon also entails (\ref{12.2}), because (\ref{12.3}) along with the inequality $\iota>\frac{1}{2}$ particularly
  ensures that $\sup_{t\in (0,\tm)} \|u(\cdot,t)\|_{L^\infty(\Om)}$ is finite. 
  Once more relying on Lemma~\ref{lem11}, going back to the first equation in (\ref{0}) and invoking the outcome
  of a Moser type iterative reasoning (\cite[Lemma A.1]{taowin_subcrit}) we finally obtain also (\ref{12.1}).
\qed
\proofc of Theorem~\ref{theo13}. \quad
  It is sufficient to combine Lemma~\ref{lem12} with Lemma~\ref{lem_loc}.
\qed

\bigskip

{\bf Acknowledgment.} \quad
The second author acknowledges support of the {\em Deutsche Forschungsgemeinschaft} (Project No.~462888149).\\

{\bf Conflict of interest statement.} \quad
On behalf of all authors, the corresponding author states that there is no conflict of interest.\\

{\bf Data availability statement.} \quad
No new data were created or analysed during this study. Data sharing is not applicable to this article.


\end{document}